\documentclass{amsart}
\usepackage{amsmath}
\usepackage{amsfonts}
\usepackage{amssymb}

\newtheorem{Theorem}{Theorem}[section]

\newtheorem{Proposition}[Theorem]{Proposition}

\newtheorem{Definition}[Theorem]{Definition}

\numberwithin{equation}{section}

\newcommand{\R}{\mathbb R}

\newcommand{\N}{\mathbb N}

\newcommand{\C}{\mathcal{C}}
\newcommand{\F}{\mathcal{F}}

\author{Jean-Pierre Magnot}
\title{Differentiation on spaces of triangulations and optimized triangulations. }
\address{Lyc\'ee Jeanne d'Arc \\ Avenue de Grande Bretagne \\ F-63000 Clermont-Ferrand}
\email{jean-pierr.magnot@ac-clermont.fr}
\begin{document}
	\begin{abstract}
		We describe a smooth structure, called Fr\"olicher space, on CW complexes and spaces of triangulations. This structure enables differential methods for e.g. minimization of functionnals. As an application, we exhibit how an optimized triangulation can be obtained 
		in order to solve a standard PDE.
	\end{abstract}
	\maketitle
	\section*{Introduction}
	\section{Fr\"olicher structures on CW complexes}
	\subsection{Fr\"olicher spaces}
	\label{1.1}

	\begin{Definition} $\bullet$ A \textbf{Fr\"olicher} space is a triple
		$(X,\F,\C)$ such that
		
		- $\C$ is a set of paths $\R\rightarrow X$,
		
		- A function $f:X\rightarrow\R$ is in $\F$ if and only if for any
		$c\in\C$, $f\circ c\in C^{\infty}(\R,\R)$;
		
		- A path $c:\R\rightarrow X$ is in $\C$ (i.e. is a \textbf{contour})
		if and only if for any $f\in\F$, $f\circ c\in C^{\infty}(\R,\R)$.
		
		\vskip 5pt $\bullet$ Let $(X,\F,\C)$ et $(X',\F',\C')$ be two
		Fr\"olicher spaces, a map $f:X\rightarrow X'$ is \textbf{differentiable}
		(=smooth) if and only if $\F'\circ f\circ\C\in C^{\infty}(\R,\R)$.
	\end{Definition}
	
	Any family of maps $\F_{g}$ from $X$ to $\R$ generate a Fr\"olicher
	structure $(X,\F,\C)$, setting \cite{KM,Ma2015}:
	
	- $\C(\F_g)=\{c:\R\rightarrow X\hbox{ such that }\F_{g}\circ c\subset C^{\infty}(\R,\R)\}$
	
	- $\F(\C(\F_g))=\{f:X\rightarrow\R\hbox{ such that }f\circ\C\subset C^{\infty}(\R,\R)\}.$
	
	A Fr\"olicher space carries a natural topology,
	which is the pull-back topology of $\R$ via $\F$, see e.g. \cite{BT2014}. In the case of
	a finite dimensional differentiable manifold, the underlying topology
	of the Fr\"olicher structure is the same as the manifold topology. In
	the infinite dimensional case, these two topologies differ very often.

	\vskip 6pt
	The classical properties of a tangent space are not automatically checked in the case of Fr\"olicher spaces following \cite{DN2007-1}.  There are two possible tangent spaces:
	
	-the \textbf{internal} tangent space, defined by the derivatives of smooth paths,
	
	- the \textbf{external} tangent space, made of derivations on $\R-$valued smooth maps. 
	
	The internal tangent space is not necessarily a vector space, were as the external tangent space is. For example, consider the diffeological subspace of $\R^2$ made of the two lines $y=x$ and $y = -x.$ They cross at the origin, the external tangent space is of dimension 2, the internal tangent space is made of two directions that cannot be combined by addition. This is why some authors sometimes call the internal tangent space by ``tangent cone''.

	\begin{Proposition} \label{quotient}
		Let $X$ be a Fr\"olicher space equipped with an equivalence relation $\mathcal{R}.$
		The quotient map $\pi : X \rightarrow X/\mathcal{R}$ defines by push-forward a diffeology on $X /\mathcal{R}.$
	\end{Proposition}

	We can even state the same results in the case of infinite products,
	in a very trivial way by taking the cartesian products 
	 of the contours.  Let us now give the description of what happens
	for projective limits of Fr\"olicher spaces.

	\subsection{Topological gluing of Fr\"olicher spaces}\label{glue}
	Let us assume that $X$ is a topological space, and that there is a collection $\{(X_i,\F_i,\C_i)\}_{i \in I}$ of Fr\"olicher spaces, 
	together with continuous maps $f_i: X_i \rightarrow X.$
	Then we can define a Fr\"olicher structures on $X$
	setting $$\F_{I,0} = \{f \in C^0(X,\R) | \forall i \in I, \quad f \circ f_i \circ \C_i \subset C^\infty(\R,\R)\},$$ wa define $\C_I$ the contours generated by the family $\F_{I,0},$ and  $\F_I = \F(\C_I).$
	
	\vskip 12pt
	\noindent
	\textbf{Example: Fr\"olicher structure on CW complexes.}
	
	A CW-complex is a topological space built by induction:
	
	- \underline{1st step:}  gluing a family of intervals  of the type $[0;1]$ along their broder, in order to have a graph. We get a $1-$CW complex.
	
	-\underline{Assume that we have obtained a $n-$CW complex $X_n$} We glue topologicaly some $(n+1)$ disks $$ D_{n+1} = \{ x \in \R^{n+1} | ||x||\leq 1 \}$$ identifying the plots of the  
	border $S^n$ with some plots of $X_n.$
	
	\vskip 6pt
	Since at each step we glue together some Fr\"olicher space, the obtained CW complex is a Fr\"olicher space.
	
	\subsection{Fr\"olicher structure of a triangulation}
	
	Let $M$ be a smooth manifold for dimension $n.$ Let $$\Delta_n= \{(x_0,...,x_n)\in \R_+^{n+1} | x_0 +... +x_n = 1\}$$ be the $n-$ simplex. a triangulation of $M$ is a family $\sigma = (\sigma_i)_{i \in I}$
	where $I \subset \N$ is a set of indexes, finite or infinite, each $\sigma_i$ is a smooth map $\Delta_n \rightarrow M,$ and such that:
	\begin{enumerate}
		\item $\forall i \in I, \sigma_i$ is an embedding.
		\item $\bigcup_{i \in I }\sigma_i(\Delta_n) = M.$ (open covering)
		\item $\forall (i,j) \in I^2,$ $\sigma_i(\Delta_n) \cap \sigma_j(\Delta_n) \subset  \sigma_i(\partial \Delta_n) \cap \sigma_j(\partial \Delta_n).$ (intersection along the borders)
		\item  $\forall (i,j) \in I^2$ such that $ D_{i,j}=\sigma_i(\Delta_n) \cap \sigma_j(\Delta_n) \neq \emptyset,$ for each $(n-1)$-face $F$ of $D_{i,j},$ the ``transition maps" $\sigma_j^{-1} \circ \sigma_i : \sigma_i^{-1}(F) \rightarrow \sigma_j^{-1}(F)$ are affine maps.  
	\end{enumerate}  
	Under these (well-known) conditions, we can apply section \ref{glue}, and equip the triangulated manifold $(M,\sigma)$ with a Fr\"olicher structure $(\F,\C),$ generated by the smooth maps $\sigma_i.$ The following result is obtained from the construction of $\F$ and $\C:$
	\begin{Theorem}
		The inclusion $(M,\F,\C)\rightarrow M$ is smooth.
	\end{Theorem}
	\section{Smooth structures on the space of triangulations}
	
	\subsection{The Fr\"olicher structure}
	
	Now, let us fix the set of indexes $I$ and fix a so-called \textbf{model triangulation} $\sigma.$ We note by $\mathcal{T}_\sigma$ the set of triangulations $\sigma'$ of $M$ such that the domains $ D_{i,j}=\sigma_i(\Delta_n) \cap \sigma_j(\Delta_n)$ and $ D'_{i,j}=\sigma'_i(\Delta_n) \cap \sigma'_j(\Delta_n)$ are diffeomorphic (in the Fr\"olicher category). Based on the properties described in \cite{KM}, since $\mathcal{T}_\sigma \subset C^\infty(\Delta_n, M)^I,$ we can equip  $\mathcal{T}_\sigma$ with the trace Fr\"olicher structure, in other words, the Fr\"olicher structure on $\mathcal{T}_\sigma$ whose generating family of contours $\C$ are the contours in $C^\infty(\Delta_n, M)^I$ which lie in $\mathcal{T}_\sigma.$ In practice, as we shall see below, we often consider a subspace of $\mathcal{T}_\sigma.$
	
	\subsection{Application: optimized approximations}
	We wish here to describe applications to optimization of $H^1-$ approximation of functions and differential forms, in the sense of \cite{Wh}, which is commonly used in numerical analysis. $M$ is assumed Riemannian. In the two examples, we consider the Sobolev space $H^1(M,\R)$ or $H^1(\Omega^*(M,\R))$, both equipped with self-adjoint positive injective pseudo-differential operator $A$ of order 2. With this setting, it is well-known that $(A.,.)_{L^2}$ is a $H^1-$scalar product, where $(.,.)_{L^2}$ is the $L^2-$scalar product. For example, we can assume that $A-Id$ is the Laplace-Beltrami operator on $H^1(M,\R),$ or the Hodge Laplacian on $H^1(\Omega^*(M,\R)).$ These approximations are the base of the finite elements method.
	\subsubsection{Approximation of functions in $H^1(M,\R)$}
	For a fixed map $f \in C^0(M,\R),$ the approximation $f_{\sigma'}$ of $f$ with respect to the triangulation $\sigma' \in \mathcal{T}_\sigma$ is given by the following constraints:
	\begin{itemize}
		\item $f_{\sigma'}= f$ at the 0-vertex of the triangulation
	\item $\forall \sigma'_i \in \sigma', f_{\sigma'} \circ \sigma_i$ is an affine map $\Delta_n \rightarrow \R.$
	\end{itemize}
	
	\begin{Definition}
		An optimized $\mathcal{T}_\sigma$ triangulation $\sigma_{op} (f)$ is a triangulation that minimize $$\Phi: \sigma' \mapsto (A(f_{\sigma'} - f),f_{\sigma'} - f )_{L^2}.$$ 
	\end{Definition}
	
	This is a smooth map, and the condition $D\Phi = 0$ is necessary at the minimizing point $\sigma_{op}.$
	\subsubsection{Approximation of differential forms in $H^1(\Omega^*(M,\R))$} 
	Let $\alpha \in \Omega^k(M,\R).$ Following \cite{Wh}, for any triangulation $\sigma' \in \mathcal{T}_\sigma,$  there exists a discretization of $\alpha$ with respect to $\sigma'$ which can be realized as a $H^1$-differential form.
	Then, one can get, the same way, an optimized triangulation $\sigma_{op}.$  
	
	\section{Optimized triangulations in $H^1$ approximation}
	 If one considers the whole space $\mathcal{T}_\sigma,$ one quicly see that the minimization functionnal furnished only a reparametrization of the function $f.$
	 
	 Let us consider the following example: let $f(x)=\frac{x + x^2}{2}$ defined on $[0;1],$ and let us fix the model triangulation as the trivial parametrization $[0;1] \rightarrow [0;1].$ An optimal triangulation $\sigma_{op} : [0;1] \rightarrow [0;1]$ is then a smooth map such that $$\sigma_{op} = f(x).$$ This shows that, if we wish one minimize $\Phi$ on  $\mathcal{T}_\sigma$,
	 
	 - we have got here an optimized triangulation $\sigma_{op} = f(x),$ which is a toy example of what can happen when, for example, the first derivative of $f$ is Lipschitz and when $M$ is a compact domain;  
	 
	 - defining an optimal triangulation leads to choose a function $\sigma_{op}$ in a too large class of functions.
	 
	 \vskip 12pt
	 In the finite elements method, finding a base triangulation is not the main feature, and one needs to get quickly a ``not so bad'' triangulation. This is why one usually considers, when $M$ is a bounded domain of an Euclidian space, the space of affine triangulations: 
	  $$Aff\mathcal{T}_\sigma = \left\{ \sigma' \in \mathcal{T}_\sigma | \forall i , \sigma_i' \hbox{ is affine } \right\}.$$
	  
	  \section{Conclusion and perspectives}
	  Under the lights of this approach, which was not fully justified till the definition of the Fr\"olicher structure of a triangulation in \cite{N2009} and the remark that the space of triangulations could be equipped with a smooth structure in the present communication, we can sketch two promising directions:
	  
	  - the systemetic review of existing procedures for choosing a triangulation, which can be very sophisticated, but often justified filling technical gaps by intuition. The present setting may furnish a technical tool to fully justify existing methods.
	  
	  - the exploration of the variation of functionnal integrals in the space of connections, which are often defined via cylindrical integrals over $H^1$ approximations with respect to a sequence of triangulations.

\end{document}